\documentclass[12pt]{amsart}
\usepackage{amsmath,amstext,amsthm,amssymb}

\numberwithin{equation}{section}

\newcommand{\field}[1]{\mathbb{#1}} 
 
\newcommand{\R}{\field{R}} 
\newcommand{\C}{\field{C}} 
 
\newcommand{\caligra}[1]{\mathcal{#1}} 
\newcommand{\ci}{\caligra{I}} 
\newcommand{\co}{\caligra{O}} 
\newcommand{\cC}{\caligra{C}} 
\newcommand{\cu}{\caligra{U}} 
 
\newcommand{\f}{\varphi} 
\newcommand{\om}{\omega} 
 
\newcommand{\lam}{\lambda} 
\newcommand{\g}{\gamma} 
\newcommand{\G}{\Gamma^{\ve}} 
\newcommand{\ve}{\varepsilon} 
\newcommand\db{\bar\partial}
\newcommand\dbad{\db^{\ast}} 
\newcommand\pa{\partial}  
\newcommand{\curv}[1]{\imath\mathbf{c}{#1}} 
\newcommand\ce{\curv(E)} 
\newcommand\mef{M_{E}(\f)} 
\newcommand\memf{M_{E}(-\f)}  
\newcommand\mome{M^{\prime}_{\om}(E)} 
\newcommand\norm[1]{\left\|{#1}\right\|} 
\newcommand\abs[1]{\vert{#1}\vert}

\newcommand{\latk}{\tfrac{1}{k}\Delta^{\prime\prime}_{k,\varepsilon}}
\newcommand{\latd}{\tfrac{1}{k}\Delta^{\prime\prime}_{k,\varepsilon} 
\upharpoonright_{D(\ve/2)}}
\newcommand{\nuk}{N^j\left(\lam,\latk\right)}
\newcommand{\nud}{N^j\left(\mu,\latd\right)} 
 
\DeclareMathOperator{\degtr}{deg\,tr} 
\DeclareMathOperator{\codim}{codim} 
\DeclareMathOperator{\Dom}{Dom} 
\DeclareMathOperator{\Ran}{Ran} 
 
\newtheorem{theorem}{{\sc Theorem}}[section]
\newtheorem{lemma}[theorem]{{\sc Lemma}}
\newtheorem{cor}[theorem]{{\sc Corollary}}
\newtheorem{prop}[theorem]{{\sc Proposition}}
\newtheorem{property}[theorem]{{\sc Property}} 
\newtheorem*{defo}{{\sc Stability Theorem}} 
\newtheorem*{conc}{{\sc Existence Criterion}} 

\theoremstyle{remark}

\newcommand{\comment}[1]{}

\begin{document}
\title[Existence of holomorphic sections and perturbation]{Existence of holomorphic sections and 
perturbation of positive line bundles over $q$--concave manifolds}  
%\title{On the deformation of certain pseudoconcave complex manifolds}
\author{George Marinescu}
\date{February 5, 2002}
\address{Institut f\"ur Mathemathik, 
Humboldt--Universit\"at zu Berlin,
Unter den Linden 6,
10099 Berlin, Deutschland}

\email{george@\,mathematik.hu-berlin.de}

\address{Institute of Mathematics of the Romanian Academy, 
PO Box 1--764, RO--70700,
Bucharest, Romania}

\maketitle 

Let $X$ be a compact complex manifold with strongly pseudoconcave
boundary. The question of projectively embedding $X$
has been studied in \cite{An}, \cite{AT}, \cite{AS}, \cite{EH}, 
\cite{Gri}. In particular the generalization of Kodaira's embedding theorem
would give an intrinsic characterisation of projective pseudoconcave manifolds
in terms of positive line bundles. 
Our model is the case when $X$ admits a positively embedded 
(i.e. with positive normal bundle) smooth compact divisor $Z$.
By a rigidity theorem of Griffiths \cite{Gri} 
we infer that global sections in high tensor powers of the 
the associated bundle $[Z]$ embed a small neighbourhood of $Z$
in the projective space. In particular
$X$ has a maximal number of independent meromorphic functions.

We will be concerned in the sequel with general positive line bundles.
Since for $\dim X=2$ there exists 
examples of strongly pseudoconcave manifolds which possess positive line 
bundles but cannot be even compactified (see Andreotti--Siu\cite{AS}) we will restrict 
ourselves to the case $\dim X\geqslant 3$. For an analysis of the case $\dim X=2$, see
Epstein--Henkin \cite{EH}.

A first step towards the 
Kodaira embedding of $X$ is to find holomorphic sections in the tensor
powers $E^k$ of a positive line bundle $E$ over $X$. We will prove an
existence criterion giving a lower bound for $\dim H^0(X,E^k)$ in terms
of geometric data such as the Levi form of $\pa X$ and the curvature of $E$.
As a corollary we see that, roughly speaking, if the volume of $X$ in the 
metric $\ce$ exceeds the volume of $\pa X$ times a constant expressing the 
size of the Levi form and of the curvature $\ce$ near the boundary,
%calculated in the same metric
the ring $\oplus H^0(X,E^k)$ contains local coordinates for each point
outside a proper analytic set of $X$.
Since for $\dim X\geqslant 3$ there exists a compactification $\hat{X}$
which contains $X$ as an open set we deduce that
$\hat{X}$ is Moishezon. We remind that by an important theorem of Moishezon
a compact complex manifold with a maximal number of independent meromorphic
functions (this being the definition of Moishezon manifolds) is not far 
from being projective. 
Namely, there exists a proper projective modification
of $\hat{X}$. This implies that modifiying $X$ along a proper analytic set
(which may cut the boundary) we obtain an open set in a projective manifold.   
Note that the main result of Andreotti--Tomassini \cite{AT} (see also \cite{AS}) 
says that if $\oplus H^0(X,E^k)$ separates points and gives local coordinates on $X$
there exists a projective manifold $\hat{X}$ containing $X$ as an open set.
Our aim is to find geometric conditions which imply the hypothesis in
Andreotti--Tomassini theorem.

In the pesent paper we work actually with a more general class of manifolds,
namely $q$--concave manifolds. In this terminology, 
strongly pseudoconcave manifolds correspond to $1$--concave manifolds.
As application of the existence theorem we prove a stability property for
certain $q$--concave manifolds.
Let us consider the complement $X$ of a sufficiently small neighbourhood
of a submanifold of codimension $\geqslant 3$ in a projective manifold.
Assume that we perform a small perturbation of the complex stucture of $X$ 
such that along a (not necessaraly compact) smooth divisor the structure 
remains unchanged. Then the resulting manifold still
has a maximal number of meromorphic functions. 
If moreover the canonical bundle is positive, any small perturbation
suffices for the result to hold.

\section*{Aknowledgements}
A version of the existence criterion for holomorphic sections was obtained in a 
joint manuscript with G. Henkin, using the original method of Siu. The constants
were however not so explicit as in the present paper.

\comment{
We show that for small perturbation of the complex structure of the
complement of a small neighbourhood of a submanifold of codimension 
$\geqslant 3$ in a projective manifold such that along a 
(not necesseraly compact) smooth divisor the structure remains unchanged
the resulting manifold still
has a maximal number of meromorphic functions.

We show that for small perturbation of the complex structure of the 
complement of a small neighbourhood of a submanifold of codimension 
$\geqslant 3$ in a projective manifold of general type the resulting manifold still
has a maximal number of meromorphic functions. 
}

%%%%%%%%%%%%%%%%%%%%%%%%%%%%%%%%%%%%%%%%%%%%%%%%%%%%%%%%%%%%%%%%%% 
\section{Description of the results} 
In this paper we shall be concerned with %the deformation of  
{\it pseudoconcave} (for short {\it concave}) complex manifolds. 
We understand concavity in the sense of Andreotti--Grauert \cite{AG}. 
A manifold $X$ of dimension $n$ is called $q$--\,concave if there exists  
a smooth function $\f:X\longrightarrow (a,b\,]$ where  
$a=\inf\f\in\{-\infty\}\cup\R$, $b\in\R$, such that $X_c:=\{\f>c\}\Subset X$ for  
all $c\in (a,b\,]$ and $\imath\pa\db\f$ has at least $n-q+1$ positive  
eigenvalues outside an exceptional compact set $K$. 
The prime examples of such manifolds arise as complements of complex  
submanifolds of compact manifolds. More precisely, let $M$ be a compact 
complex manifold and $A\subset M$ of dimension $q$. Then $M\smallsetminus A$ 
is $(q+1)$--concave (see \S 4.). 
It is well known (see \cite{An}) that for  
a $q$--concave manifold $X$ ($q\leqslant n-1$) the transcedence degree  
$\degtr\mathcal K(X)$ of 
the meromorphic function field is at most the complex dimension of $X$. 
In analogy to the corresponding notion for compact manifolds we say that 
a $q$--concave manifold is {\it Moishezon} if $\degtr\mathcal K(X)= 
\dim_{\C}X$. 
 
Let us consider now a projective manifold $M$, a submanifold $A\subset M$ 
and the concave manifold $X:=M\smallsetminus A$. Our aim is to study to what  
extent small deformations of the sublevel sets $X_c$ for small values  
of $c>\inf\f$ (i.e. for $X_c$ close to $X$) give rise to concave Moishezon  
manifolds.  
%Why only Moishezon? 
As a matter of fact we may consider 
small neighbourhoods $V$ of $A$, which means that $X_c\subset 
M\smallsetminus V$ for small $c>\inf\f$. Then  $M\smallsetminus V$ is  
pseudoconcave in the sense of Andreotti and the notion of Moishezon 
manifold still makes sense (see \cite{An}). 
 
\begin{defo}\label{defo} 
Let $M$ be a compact projective manifold and let $Z$ be an ample smooth divisor.
Let $A\subset M$ be a complex submanifold of codimension at least $3$.
Then for any sufficiently small neighbourhood $V$ of $A$ and for any 
sufficiently small deformation of the complex structure of 
$M\smallsetminus V$ leaving $T(Z)$ invariant, the manifold $M\smallsetminus V$ with the new 
structure is a pseudoconcave Moishezon manifold. If the canonical bundle $K_M$ is positive,
the statement holds for any small enough perturbation.
\end{defo}

Let us note that the sublevel sets $X_c$ are also $q$--concave manifolds 
if $K\Subset X_c$. 
Our method is based on $L^2$ estimates for $(0,1)$--forms on $X_c$  
and for this reason we need at least $3$ positive eigenvalues for  
$\imath\pa\db\f$, that is $\codim A=n-q\geqslant 3$. 
While this condition may seem technical we can explain it as follows. 
The existence of $L^2$ estimates for $(0,1)$--forms imply the  
finiteness of the first cohomology group $H^1(X,F)$ for holomorphic  
vector bundles $F$ over $X$. By the Andreotti--Grauert theory 
we have $\dim H^p(X,F)<\infty$ for $p\leqslant n-(q+1)-1=n-q-2$ 
and $\dim H^p(X,F)=\infty$ for $p=n-q-1$. Therefore we have to impose 
$n-q-1>1$ i.e. $n-q>2$. 
 
An immediate consequence is the following. 
\begin{cor} 
Let $M$ be a compact projective manifold and let let $Z$ be an ample smooth divisor.
Let $A\subset M$ be a complex submanifold of codimension at least $3$.
Then for any sufficiently small neighbourhood $V$ of $A$ and for any 
deformation of the complex structure of $M$ which 
is sufficiently small on $M\smallsetminus V$ and leaves $T(Z)$ invariant,
the manifold $M$ with the new structure is Moishezon.
\end{cor} 
In order to prove the Stability Theorem we need a differential geometric 
criterion for a $q$--concave manifold to be Moishezon. For compact 
manifolds the type of results we need were proved by Siu \cite{Si} 
and Demailly \cite{De}. They derive asymptotic Morse inequalities 
for the cohomology groups with values in the tensor powers of a  
holomorphic line bundle. For non--compact manifolds the Morse 
inequalities were used by Nadel--Tsuji \cite{NT} to prove the 
quasi--projectivity 
of very strongly $(n-2)$--concave manifolds of dimension $n$ 
which possess a complete K\"ahler metric with 
$\operatorname{Ric}\om<0$ and whose universal covering is Stein. 
In \cite{Ma}, \cite{TCM} we considered Morse inequalities for general 
$1$--concave manifolds with application to the deformation of the  
complex structure of compact complex spaces with isolated singularities. 
In the sequel we study $q$--concave manifolds and give an estimate from  
below of the dimension of the space of holomorphic sections with values 
in a positive line bundle (see the Existence Criterion below). 
An important feature of our estimate is the presence of a negative 
boundary term which expresses the obstruction to finding holomorphic  
sections. 
 
We need some preparations and notations in order to state the result. 
Let $X$ be a $q$--concave manifold with exhaustion function $\f$. 
If $\pa X_c$ is smooth the Levi form of $\pa X_c$ 
has at least $n-q-1$ negative eigenvalues (since the defining function  
for $X_c$ is $c-\f$). Therefore the following setting may be considered. 
 
Let $D\Subset X$ be a smooth domain in a complex manifold $X$ such that 
the Levi form of $\pa D$ has at least $2$ negative eigenvalues. 
Then we can choose a defining function $\f$ for $D$ which is smooth on 
$\overline{D}$, $D=\{\f<0\}$ and $\pa\db\f$ has at least $3$ negative  
eigenvalues. We can in fact modify a defining function in order to get 
an extra negative eigenvalue in the complex normal direction to $\pa D$. 
In the following we keep the function $\f$ fixed.  
 
We introduce a hermitian metric $\om=\om_{\f}$ in the neighbourhood of  
$\overline{D}$ such that in a neighbourhood $V$ of $\pa D$  
the following property holds: 
\begin{property}\label{metric} 
The first $3$ eigenvalues of $\imath\pa\db\f$ with respect  
to $\om$ are at most $-2n+3$  
and all others are at most $1$. 
\end{property} 
%We need the following general notation.  

Finally set $dS_E$ for the volume form of $\pa D$ in the induced metric 
from $\ce$ and $|d\f|_E$ for the norm of $d\f$ in the metric associated to 
$\ce$. 

We can state the estimate for the dimension of the space holomorphic  
sections on the concave domain $D$. 
\begin{conc}\label{conc} 
Let $D\Subset X$ be a smooth domain in a complex manifold $X$ such that the 
Levi form of $\pa D$ possesses at least $2$ negative eigenvalues. 
Let $E$ be a holomorphic line bundle on $X$ which is assumed to be positive  
on a neighbourhood of $\overline{D}$. Then 
\begin{equation}\label{conc1} 
\liminf_{k\longrightarrow\infty}\;k^{-n}\dim H^0(D,E^k)\geqslant 
\int_D\left(\tfrac{\imath}{2\pi}\mathbf{c}(E)\right)^n 
-C(\f,E)\int_{\pa D}\frac{dS_{E}}{|d\f|_E} 
\end{equation} 
The constant $C(\f,E)$ depends explicitely on the curvature of $E$ and on the  
Levi form $\imath\pa\db\f$ {\rm(}cf. \eqref{const}{\rm)}. 
\end{conc} 

%%%%%%%%%%%%%%%%%%%%%%%%%%%%%%%%%%%%%%%%%%%%%%%%%%%%%%%%%%%%%%%%%%%%%%%%%%%%%%%%%%5
\section{Proof of the Existence Criterion}
A familiar method of producing holomorphic sections in a positive bundle $E$ 
is the use of $L^2$ estimates for $\db$ of Andreotti--Vesentini and H\"ormander 
(see e.g. \cite{AV}, \cite{De1} and \cite{Oh}). The $L^2$ estimates may be established 
equally in the case of pseudoconvex and pseudoconcave manifolds 
by introducing `weights' (i.e. changing the hermitian metric on the bundle) 
that reflect the convexity or concavity of the manifold. 
The problem is that, for pseudoconcave manifolds the positivity 
is lost by this procedure (in contrast to the pseudoconvex case). 
There is however a strategy of finding holomorphic sections in non--positive 
hermitian bundles which has been introduced by Siu \cite{Si} for semipositive 
line bundles and then generalized 
by Demailly \cite{De} into his asymptotic Morse inequalities. 
The main ingredient is a Weyl type formula describing the semiclassical 
behaviour of the $\db$--laplacian on the tensor powers $E^k$. 
The first applications for non--compact manifolds appear in Nadel--Tsuji 
\cite{NT} and Bouche \cite{Bou}. 
  
We proceed as follows. In a first instance we find a good $L^2$--estimate 
for the $(0,1)$--forms with values in $E^k$. Then following Bouche \cite{Bou} 
we compare the spectrum of the Laplace operator on $D$ (for a complete metric) 
with the spectrum of the Dirichlet problem over a smaller domain $D(\ve/2)$ 
which is a set of points of $D$ at distance less than $\sqrt{\ve/2\,}$ 
times a certain constant 
from $\pa D$ (see \eqref{not} for the precise definition). 
On $D(\ve/2)$ we can use Demailly's spectral formula and get a lower bound for  
the dimension of the space of sections in $E^k$ for large $k$. 
We shall need the full strength of Demailly's result since the curvature 
of the changed metric has negative eigenvalues. 
In the last step we apply the results to metrics 
which approximate the positive metric on $E$ in the interior of the manifold. 
In the process of approximation the set where the curvature has a negative  
part concentrates to the boundary $\pa D$ and is responsible for the negative boundary 
term in the final estimate of the Existence Criterion. 

We begin by setting some notations and defining the constant $C(\varphi,E)$.
  
Let $\eta$ a hermitian metric 
on $X$, $\Phi$ a real $(1,1)$--form and $K$ a compact set in $X$. 
We set: 
$$M_{\eta}(\Phi,K)=\sup_{x\in K}\;\sup_{v\in T_{x}X\smallsetminus\{0\}} 
\frac{\Phi(v,v)}{\eta(v,v)}\,,$$ 
the supremum over $K$ of the highest eigenvalue of $\Phi$ with respect to 
$\eta$. 
In hindsight to our previous situation denote: 
\begin{align*} 
&M_{E}(\f)=M_{\ce}(\imath\pa\db\f,\overline{D})\\ 
&M_{E}(-\f)=M_{\ce}(-\imath\pa\db\f,\overline{D})\\ 
&M_{\om}(E)=M_{\om}(\ce,\overline{D})\\ 
&M^{\prime}_{\om}(E)=1+2(n-1)\,M_{\om}(\ce,\overline{D})\\ 
&M_{E}(\pa\f)=M_{\ce}(\imath\pa\f\wedge\db\f,\pa{D}) 
\end{align*} 
which represent the relative size of the respective $(1,1)$--forms. 
We also put: 
\begin{align*} 
C_1&=\sqrt{2\memf\mome}-1\\ 
C_2&=2\memf\mome-1\\ 
C_3&=2\mef\mome+1\\ 
C_4&=2\mome M_{E}(\pa\f) 
\end{align*} 

The definition of $C(\f,E)$ is then
\begin{equation}\label{const} 
C(\f,E)=(2\pi)^{-n}\,C_1\,C_2\,C_3^{n-2}\,C_4\,. 
\end{equation} 
 
Let $\g_1\leqslant\gamma_2\leqslant\dotsb\leqslant\gamma_n$ be the eigenvalues  
of $\imath\pa\db\f$ with respect to $\om$. We have chosen $\om$ such that 
(see Property \ref{metric}) in a neighbourhood $V$ of $\pa D$,  
\begin{subequations} 
\begin{align} 
\g_1 &\leqslant\g_2\leqslant\g_3\leqslant -2n+3 \,, \label{hip1}\\ 
\g_n &\leqslant 1\,.\label{hip2}  
\end{align} 
\end{subequations} 
Let $\chi:(-\infty,0)\longrightarrow\R$, $\chi(t)=t^{-2}$. We consider 
the complete metric: 
\begin{equation} 
\om_0=\om+\chi(\f)\pa\f\wedge\db\f 
\end{equation} 
which grows as $\f^{-2}$ in the normal direction to $\pa D$. Along the 
fibers of $E$ we introduce the metric: 
\begin{equation} 
h_\ve=h\,\exp\left(-\ve\int_{\inf\f}^{\f}\chi(t)\,dt\right) 
\end{equation} 
where $h$ is the given metric on $E$ (for which $\ce$ is positive). 
The curvature of $h_\ve$ is 
$$ 
\curv(E,h_\ve)=\ce+\imath\ve\chi(\f)\pa\db\f+ 
\imath\ve\chi^{\prime}(\f)\pa\f\wedge\db\f 
$$ 
We evaluate the eigenvalues of $\curv(E,h_\ve)$ with respect to  
$\om_0$ with the goal to apply the Bochner--Kodaira formula. Denote by 
$\g^0_1\leqslant\g^0_2\leqslant\dotsb\leqslant\g^0_n$ the eigenvalues 
of $\imath\pa\db\f$ and 
$\G_1\leqslant\G_2\leqslant\dotsb\leqslant\G_n$ the eigenvalues of  
$\imath\ve\chi(\f)\pa\db\f+\imath\ve\chi^{\prime}(\f)\pa\f\wedge\db\f$  
with respect to $\om_0$. The minimum--maximum principle yields 
\begin{subequations} 
\begin{align} 
\g_1&\leqslant\g^0_1\leqslant\g_2\leqslant\g^0_2\leqslant\g_3\leqslant  
-2n+3 &&\quad\text{by \eqref{hip1}\,,}\label{eqA}\\ 
\g^0_3&<0 &&\quad\text{since $\g_3<0$\,,}\label{eqB}\\ 
\g^0_j &\leqslant\max\{\g_n,0\}\leqslant 1  
\quad\text{for $4\leqslant j\leqslant n$\,,}\label{eqC} 
&&\quad\text{by \eqref{hip2}\,.} 
\end{align} 
\end{subequations} 
on $V$. 
It is also easy to see that the highest eigenvalue of  
$\imath\chi^{\prime}(\f)\pa\f\wedge\db\f$ with respect to $\om_0$ satisfies 
\begin{equation}\label{eqD} 
\sup_{v\in T_{x}X\smallsetminus\{0\}} 
\frac{\imath\chi^{\prime}(\f)\pa\f\wedge\db\f(v,v)}{\om_{0}(v,v)}\leqslant 
\chi(\f)\,,\quad\text{for all $x\in D$\,.} 
\end{equation} 
By \eqref{eqD} we have 
$$\G_j\leqslant\ve\chi(\f)(\g^0_j+1)$$ 
and therefore,  
\begin{align*} 
\G_1&\leqslant\G_2\leqslant (-2n+4)\ve\chi(\f)&&\quad\text{by \eqref{eqA}\,,}\\ 
\G_3&\leqslant\ve\chi(\f)&&\quad\text{by \eqref{eqB}\,,}\\ 
\G_j&\leqslant2\ve\chi(\f)\quad\text{for $4\leqslant j\leqslant n$\,,} 
&&\quad\text{by \eqref{eqC}\,.} 
\end{align*} 
Summing up we obtain 
\begin{equation}\label{sum} 
\G_2+\dotsb+\G_n\leqslant-\ve\chi(\f)\,. 
\end{equation} 
This sum will appear in the Bochner--Kodaira formula and carries the  
information about the concavity of $D$. 
 
We also have to estimate the eigenvalues of $\ce$ with respect to $\om_0$. 
We denote by $\alpha_1\leqslant\alpha_2\leqslant\dotsb\leqslant\alpha_n$ the  
eigenvalues of $\ce$ with respect to $\om$ and by 
$\alpha^0_1\leqslant\alpha^0_2\leqslant\dotsb\leqslant\alpha^0_n$ the eigenvalues of 
$\ce$ with respect to $\om_0$. It is straightforward that 
\begin{equation}\label{vp} 
\alpha^0_n\leqslant\alpha_n\leqslant M_{\om}(E)<\infty\quad\text{on $V$\,.} 
\end{equation} 
 
Since the torsion operator of $\om_0$ with respect to $\om_0$ are bounded by 
a constant $A>0$ (depending only on $\om_0$), the Bochner--Kodaira formula  
assumes the following form (see e.g. \cite{De}, \cite{Oh}): 
\begin{multline}\label{bk} 
\tfrac{3}{2}\left(\norm{\db\,u}^2+\norm{\dbad\,u}^2\right)\\ 
\geqslant\int_D\left[-k(\G_2+\dotsb+\G_n)-k(\alpha^0_2+\dotsb+\alpha^0_n)-A\chi(\f) 
\right]\abs{u}^2\,dV 
\end{multline} 
for any compactly supported $(0,1)$--form in $D$ with values in $E^k$. 
The volume form is taken with respect to $\om_0$ and the norms are  
with respect to $\om_0$ on $D$ and $h_\ve$ on $E$. 
The inequalities \eqref{sum}, \eqref{vp} and \eqref{bk} entail 
\begin{equation}\label{bk1} 
\tfrac{3}{2}\left(\norm{\db\,u}^2+\norm{\dbad\,u}^2\right)
\geqslant\int_D\left[-k(n-1)M_{\om}(E)+k\ve\chi(\f)-A\chi(\f) 
\right]\abs{u}^2\,dV
\end{equation}
for any compactly supported $(0,1)$--form in $D$ with values in $E^k$ 
and support in $V$. We use now the term $k\ve\chi(\f)$ to absorb the  
negative terms in the left--hand side of \eqref{bk1}. 
We introduce the following notation: 
\begin{equation}\label{not} 
D(\ve)=\left\lbrace x\in D\,:\,\f(x)<-\sqrt{\ve/\mome}\right\rbrace\,. 
\end{equation} 
We may assume that $V$ contains the set $\complement{D(\ve)}$  
(for $\ve$ small enough). In the  
set $\complement{D(\ve)}$ we have $\ve\chi(\f)\geqslant\mome$ and if we choose  
$k\geqslant 2A\ve^{-1}$ we get 
$$-k(n-1)M_{\om}(E)+k\ve\chi(\f)-A\chi(\f)\geqslant\frac{k}{2}$$ 
so that \eqref{bk1} yields 
\begin{equation}\label{bk2} 
3\left(\norm{\db\,u}^2+\norm{\dbad\,u}^2\right)
\geqslant k\int_D\abs{u}^2\,dV\,,\quad\operatorname{supp}u\Subset 
\complement{D(\ve)}\,,\; k\geqslant 2A\ve^{-1} 
\end{equation} 
Since the metric $\om_0$ is complete we deduce that \eqref{bk2} holds true 
for any $(0,1)$--form $u\in\Dom{\db}\cap\Dom{\dbad}$ with support 
in $\complement{D(\ve)}$ (by the density lemma 
of Andreotti--Vesentini \cite{AV}). 
 
Estimate \eqref{bk2} is crucial for our purpose. 
In Nadel--Tsuji \cite{NT} and \cite{Ma} the spectral formula of Demailly was 
used to obtain a lower bound for the dimension of the space of  
holomorphic sections of bundles over pseudoconcave manifolds. 
After having established \eqref{bk2} we should just follow the same lines.
We give the details since we need the precise output to be able to make 
$\ve\longrightarrow 0$. 
 
By following Demailly \cite{De} we reduce the problem to estimating the size
of certain spectral spaces of the $\db$--laplacian.
Let us consider the operator $\latk$ where $\Delta^{\prime\prime}_{k,\varepsilon}=
\db\db^\ast+\db^\ast\db$ is the Laplace--Beltrami operator acting on  
$(0,j)$--forms with values in $E^k$ over $D$. The 
metrics used to construct the adjoint $\db^\ast$ are $\omega_0$ and $h_\ve$.  
Let $Q_{k,\ve}$ be the quadratic form associated to $\latk$,
that is, $Q_{k,\ve}(u)=\frac{1}{k}\left(\|\db\,u\|^2+\|\db^{\ast}\,u\|^2\right)$.
We denote by $E^j\left(\lam,\latk\right)$ the spectral projectors and 
by  
\begin{align*} 
L^j\left(\lam,\latk\right)&=\operatorname{Ran}E^j\left(\lam,\latk\right)\,,\\ 
\nuk&=\dim L^j\left(\lam,\latk\right)\,, 
\end{align*}
the spectral space and the counting function for 
the spectrum of $\latk$ on $(0,j)$--forms. 

\begin{lemma}\label{lem1} For any $\lam\geqslant 0$ and $k\geqslant 0$,
$$\dim H^0(D ,E^k\otimes K_X)+N^1\left(\lam,\latk\right) 
\geqslant N^0\left(\lam,\latk\right)\,.$$
\end{lemma}
 
\begin{proof}
Since $\latk$ commutes with $\db$ it follows that the spectral projections of 
$\latk$ commute with $\db$ too, showing thus $\db L^0(\lam,\latk)\subset L^{1}(\lam,\latk)$
and therefore we have the bounded operator
$\db_\lam: L^0\left(\lam,\latk\right)\longrightarrow L^1\left(\lam,\latk\right)$
where $\db_\lam$ denotes the restriction of $\db$ (by the definition of  
$L^0\left(\lam,\latk\right)$,
$\db_\lam$ is bounded by $k\lam$).  
The assertion is a consequence of the following obvious relations: 
\begin{align*}
N^0\left(\lam,\latk\right)&=\dim\ker\db_\lam+\dim\operatorname{Ran}\db_\lam\,,\\ 
\dim\Ran\db_\lam&\leqslant N^1\left(\lam,\latk\right)\,,\\ 
\dim\ker\db_\lam&\leqslant\dim H^0(D ,E^k\otimes K_X)\,, 
\end{align*} 
where the last line follows from the fact that the kernel of  
$\db_\lam$ consits of holomorphic sections.
\end{proof} 
%Thus the inequality of Lemma \ref{lem1} assumes the form 
%all spaces appearing in the statement are finite dimensional and we 
%have the desired inequality. 
%Moreover, since $D$ is pseudoconcave,
%$\dim H^0(D ,E^k\otimes K_X)<\infty$ (see \cite{An}) so $\dim\ker\db_\lam<\infty$\,.
By the previous lemma we have to estimate $N^1(\lam,\latk)$ from above and then
$N^0\left(\lam,\latk\right)$ from below.  
%We do this thanks to a remark of Witten (see \cite{Wi},  
%p. 666 and \cite{He}, Lemma 2.1): 
%the $L^2$ norm of the eigenforms of $\latk$ on $(n,j)$--forms concentrates asymptotically 
%for $k\longrightarrow\infty$ on the critical set $D (j)$.
%In the original setting of classical Morse theory  
%the r\^ole of the curvature is played by the hessian of a Morse function $f$ 
%and the eigenforms of the modified laplacian $\Delta_t=(d_t+d^{\ast}_t)^2$ 
%where $d_t=e^{-tf}de^{tf}$ and $t>0$, concentrate near
%the critical points of $f$ as $t\longrightarrow\infty$. 
%For the complex analysis setting
%see \cite{De}, \cite{NT} and \cite{Bou}.

In the next Lemma we show that the essential spectrum 
of $\latk$ on $(0,1)$--forms does not contain the open interval $(0,1/24)$ 
and we can compare 
the counting function on this interval with the counting function of the same operator
considered with Dirichlet boundary conditions on the domain $D(\ve/2)$ (introduced 
in \eqref{not}) and denoted $\latd$. 
In particular 
$N^1\left(\lam,\latk\right)$ is finite dimensional for $\lam<1/24$.
If $E^j\left(\mu,\latd\right)$ denote the spectral projectors of $\latd$  
on $(0,j)$--forms we let 
\begin{align*} 
L^j\left(\mu,\latd\right)&=\Ran E^j(\mu,\latd)\,,\\
\nud&=\dim L^j(\mu,\latd)\,, 
\end{align*} 
be the spectral spaces and the spectrum distribution function.
For the following lemma compare \cite[Lemma 2.1]{He} and \cite[Th\'eor\`eme 2.1]{Bou}.
\begin{lemma}
For $k$ sufficiently large the operator $\latk$ on $(0,1)$--forms has discrete spectrum
in $(0,1/24)$ and  
$$N^1\left(\lam,\latk\right)\leqslant  
N^1\left(24\,\lambda +16C_\ve k^{-1},\latd\right)$$ 
for $\lam\in(0,\ve/2)$,
where $C_\ve$ is a constant independent of $k$.
\end{lemma}
\begin{proof}
Let $\rho_\ve\in{\mathcal C}^{\infty}(D)$ such that $\rho_\ve=0$
on a closed neighbourhood of $D(\ve)$ and $\rho_\ve=1$ on $\complement{D(\ve/2)}$.
Let $u\in\Dom(Q_{k,\ve})=\Dom{\db}\cap\Dom{\dbad}$ be 
a $(0,1)$--form with values in $E^k$.  
Then $\rho_\ve u$ has support in $\complement D(\ve)$ 
and for $\rho_\ve u$ we can apply \eqref{bk2}. 
We also need the following simple estimate. 
Denote $C_\ve=6\sup |d\rho_\ve|^2<\infty$. The constant depends on $\ve$ 
(which is fixed) but not on $k$.
Then
\begin{equation}\label{simple} 
Q_{k,\ve}(\rho_\ve u)\leqslant
\tfrac{3}{2}\,Q_{k,\ve}(u)+C_\ve k^{-1}\,\norm{u}^2\,. 
\end{equation} 
Using $\norm{u}^2\leqslant2\left(\norm{\rho_\ve u}^2+\norm{(1-\rho_\ve)u}^2\right)$ 
and then applying \eqref{bk2} to $\rho_\ve u$ in conjunction with \eqref{simple}
we obtain: 
\begin{equation}\label{bk3} 
\norm{u}^2\leqslant 12\, Q_{k,\ve}(u) 
+ 8\int_{D(\ve/2)}\big|(1-\rho\,)u\big|^2\,,\quad
k\geqslant\max\left\lbrace 2A\,\ve^{-1}\,,4\,C_\ve\right\rbrace 
\end{equation}
for any $u\in\Dom(Q_{k,\ve})$.  
From relation \eqref{bk3} we infer that the spectral spaces corresponding to the
lower part of the spectrum of $\latk$ on $(0,1)$--forms can be injected 
into the spectral spaces  of $\latd$ which correspond to the
Dirichlet problem on $D(\ve/2)$.
Namely, for $\lam<1/24$, the morphism  
\begin{gather*} 
L^1\left(\lam\latk\right)\longrightarrow 
L^1\left(24\lam+16C_{\ve}k^{-1},\latd\right)\,,\\
u\longmapsto E^1\left(24\lam+16C_{\ve}k^{-1},\latd\right) (1-\rho_{\ve})u 
\end{gather*}
is injective.
In order to prove the injectivity we 
choose $u\in L^1\left(\lam,\latk\right)$, $\lam<1/24$ to the effect that
$Q_{k,\ve}(u)\leqslant\lam\|u\|^2\leqslant (1/24)\|u\|^2$. Plugging this relation in 
\eqref{bk3} we get
\begin{equation}\label{bk4} 
\|u\|^2\leqslant 16\int_{D(\ve/2)}\big|(1-\rho\,)u\big|^2
\,,\quad u\in L^1\left(\lam,\latk\right)\,,\quad\lam<1/24\,. 
\end{equation}
Let us denote by $Q_{k,D(\ve/2)}$ the quadratic form of $\latd$.
Then by \eqref{simple} and \eqref{bk4},
\begin{align*}
Q_{k,D(\ve/2)}\big((1-\rho)u\big) &\leqslant
\tfrac{3}{2}\,Q_{k,\ve}(u) + C k^{-1}\,\|u\|^2\\ 
&\leqslant
\left(24\,\lam +16\,C_{\ve}k^{-1}\right)\int_{D(\ve/2)}\big|(1-\rho\,)u\big|^2\,.
\end{align*}
Thus $E^1\left(24\lam+16C_{\ve}k^{-1},\latd\right)\,(1-\rho)u=0$
entails $(1-\rho)u=0$ so that $u=0$ by \eqref{bk4}.
\end{proof}

We obtain now a lower estimate for $N^0\left(\lam,\latk\right)$.

\begin{lemma} For $\lam<1/24$ and sufficiently large $k$ the following  
relation  holds \rm{:}
$$N^0(\lam,\latk)\geqslant N^0(\lam,\latd)\,.$$
\end{lemma}
\begin{proof} 
It is straightforward to show that the $L^2$ estimate \eqref{bk2} holds also 
for $(0,0)$--forms. Therefore by repeating the proof of Lemma 2.2 we see 
that the spectrum of $\latk$ on $(0,0)$--forms is discrete in the interval 
$(0,1/24)$. We may now apply the min-max principle to the operators 
$\latk$ and $\latd$ on this interval. Since  
$\Dom (Q_{k,\ve})\supset\Dom(Q_{k,D(\ve/2)})$ the desired result follows  
immediatly. 
\end{proof}

The asymptotic behaviour of the
spectrum distribution function for the Dirichlet problem has been 
determined explicitely by Demailly \cite{De}. 
Since for $\ve$ small enough $\pa D(\ve/2)$ has measure zero we can state 
the result as follows.
\begin{prop} [Demailly]\label{dem}
There exists a function $\nu^j_{\ve}{(\mu,x)}$ depending on the 
eigenvalues of the curvature of $(E,h_\ve)$ which is bounded on compact sets of $D$ and
right continuous in $\mu$ such that for any $\mu\in\R$
\begin{equation}\label{dem1} 
\limsup_{k\longrightarrow\infty}k^{-n}\nud\leqslant\tfrac{1}{n!}\int_{D(\ve/2)}
\nu^j_{\ve}{(\mu,x)}\,dV(x)\,.
\end{equation} 
Moreover there exists an at most countable set
$\mathcal{D}_\ve\subset\R$ such that for $\mu$ outside $\mathcal{D}_\ve$  
the limit of the left--hand
side expression exists and we have equality in \eqref{dem1}. 
\end{prop}
%We do not need the explicit form of $\nu^j_{\ve}{(\lam,x)}$.
 
\noindent 
For $\lam<(1/24)$ and sufficiently large $k$ we have 
\begin{equation} 
\dim H^0(D ,E^k) \geqslant N^0\left(\lam,\latk\right)-N^1\left(\lam,\latk\right) 
%&\geqslant N^0(\lam,\latd)-N^1\left(24\,\lambda +16C_\ve k^{-1},\latd\right) 
\end{equation} 
For $\lam<(1/24)$ and $\lam$ outside $\mathcal{D}_\ve$ we apply Proposition \ref{dem} and Lemma 2.3: 
$$\lim_{k\longrightarrow\infty} k^{-n}N^0(\lam,\latk)\geqslant\tfrac{1}{n!} 
\int_{D(\ve/2)}\nu^0_{\ve}(\lam,x)\,dV(x)\,.$$  
On the other hand given $\delta>0$ we learn from Lemma 2.2 that for large $k$ 
\begin{align*} 
N^1\left(\lam,\latk\right)&\leqslant 
N^1\left(24\,\lambda +16C_\ve k^{-1},\latd\right)\\ 
&\leqslant N^1\left(24\,\lambda +\delta,\latd\right) 
\end{align*} 
hence 
$$ 
\limsup_{k\longrightarrow\infty} k^{-n}N^1(\lam,\latk)\leqslant\tfrac{1}{n!} 
\int_{D(\ve/2)}\nu^1_{\ve}(24\lam+\delta,x)\,dV(x)\,. 
$$
%\begin{eqnarray*}\limsup_{k\longrightarrow\infty}k^{-n} N^1(12\lam +C_1 k^{-1},\latd)
%&\leqslant&\limsup_{k\longrightarrow\infty}k^{-n} N^1(12\lam +\delta,\latd)\\
%&=&\displaystyle\frac{1}{n!}\int_{\Om}\nu^1_{(E,h_\ve)}{(12\lam+\delta,x)}\,dV(x)  \\
%\end{eqnarray*}
and after letting $k$ go to infinity we can also let $\delta$ go to zero.
Using these remarks we see that for all but a countable set of $\lam$ we have
$$\liminf_{k\longrightarrow\infty}k^{-n}\dim H^0(D, E^k)\geqslant
\tfrac{1}{n!}\int_{D(\ve/2)}
\big[\nu^0_{\ve}{(\lam,x)}-\nu^1_{\ve}{(24\lam,x)}\big]\,dV(x)$$
In the latter estimate we may let $\lam\longrightarrow 0$ (through values 
outside the exeptional
countable set) and this yields, by the formulas in \cite{De} for the right--hand side
\begin{equation}\label{morse1} 
\liminf_{k\longrightarrow\infty}k^{-n}\dim H^0(D, E^k)\geqslant
\tfrac{1}{n!}\int_{D(\ve/2)(\leqslant 1,h_\ve)} 
\left(\tfrac{\imath}{2\pi}{\bf c}(E,h_\ve)\right)^n 
\end{equation} 
The set $D(\ve/2)(\leqslant 1,h_\ve)$ is the set of points in $D(\ve/2)$ 
where $\curv(E)$ is non--degenerate and has at most one negative eigenvalue. 
Thus $D(\ve/2)(\leqslant 1,h_\ve)$ splits in two sets: the set $D(\ve/2)(0,h_\ve)$ 
where $\curv(E)$ is positive definite and the set $D(\ve/2)(1,h_\ve)$ where 
$\curv(E)$ is non--degenerate and has exactly one negative eigenvalue. 
The integral in \eqref{morse1} splits accordingly into one positive and 
one negative term: 
\begin{multline}\label{morse2} 
\liminf_{k\longrightarrow\infty}k^{-n}\dim H^0(D, E^k)\geqslant
\tfrac{1}{n!}\int_{D(\ve/2)(0,h_\ve)} 
\left(\tfrac{\imath}{2\pi}{\bf c}(E,h_\ve)\right)^n \\+ 
\tfrac{1}{n!}\int_{D(\ve/2)(1,h_\ve)} 
\left(\tfrac{\imath}{2\pi}{\bf c}(E,h_\ve)\right)^n 
\end{multline} 
 
Our next task is to make $\ve\longrightarrow 0$ in \eqref{morse2}. 
For $\ve\longrightarrow 0$ the metrics $h_\ve$ converges uniformly 
to the metric $h$ of positive curvature on every compact set of $D$.  
So on any compact of $D$ we recover the integral of $\curv(E)$. 
On the other hand $D(\ve/2)$ exhausts $D$ and the sets $D(\ve/2)(1,h_\ve)$ 
concentrate to the boundary $\pa D$. 
 
\noindent 
Let us fix a compact set $L\subset D$. For sufficiently small $\ve$ we have 
$L\subset D(\ve/2)$ and  
$$ 
\int_{D(\ve/2)(0,h_\ve)} 
\left(\tfrac{\imath}{2\pi}{\bf c}(E,h_\ve)\right)^n\geqslant 
\int_{L(0,h_\ve)} 
\left(\tfrac{\imath}{2\pi}{\bf c}(E,h_\ve)\right)^n 
$$ 
We have $h_\ve\longrightarrow h$ on $L$ in the $\mathcal{C}^\infty$--topology. 
Since $L(0,h)=L$ letting $\ve\longrightarrow 0$ in the previous inequality 
yields 
\begin{equation}\label{morse3} 
\liminf_{\ve\longrightarrow 0}\int_{D(\ve/2)(0,h_\ve)} 
\left(\tfrac{\imath}{2\pi}{\bf c}(E,h_\ve)\right)^n\geqslant 
\int_{L} 
\left(\tfrac{\imath}{2\pi}{\bf c}(E,h)\right)^n 
\end{equation} 
 
Let us study the more delicate second integral in \eqref{morse2}. 
For this goal we fix on $D$ the ground metric $\om_E=\curv(E)$. 
This choice will simplify our computations. 
We denote by $\lam^{\ve}_1\leqslant\lam^{\ve}_2\leqslant\dotsb\leqslant\lam^{\ve}_n$ 
the eigenvalues of $\curv(E,h_\ve)$ with respect to $\om_E$. 
Then the integral we study is 
$$ 
I_\ve=\tfrac{1}{n!}\int_{D(\ve/2)(1,h_\ve)} 
\left(\tfrac{\imath}{2\pi}{\bf c}(E,h_\ve)\right)^n 
=\tfrac{1}{(2\pi)^n}\int_{S(\ve)}\lam^{\ve}_1\,\lam^{\ve}_2\dotsm\lam^{\ve}_n\, 
\,\om^n_E/n! 
$$ 
where the integration set is 
$$ 
S(\ve):=D(\ve/2)(1,h_\ve)=\left\lbrace x\in D(\ve/2)\,:\,\lam^{\ve}_1(x)<0<\lam^{\ve}_2(x)\right\rbrace 
$$ 
We find an upper bound for $\abs{I_\ve}$ so we determine upper bounds 
for $\abs{\lam^{\ve}_1}$, $\abs{\lam^{\ve}_2}$, \dots, $\abs{\lam^{\ve}_n}$ 
on $S(\ve)$. Since $\lam^{\ve}_1$ is negative on $S(\ve)$ we have to 
obtain a lower bound for this eigenvalue. By the min-max principle 
$$ 
\lam^{\ve}_1(x)=\min_{v\in T_x D}\frac{\left[\curv(E,h)+\imath\ve\chi(\f) 
\pa\db\f+\imath\ve\chi^{\prime}(\f)\pa\f\wedge\db\f\right](v)} 
{\curv(E)(v)}\,. 
$$ 
We use now $\imath\ve\chi^{\prime}(\f)\pa\f\wedge\db\f(v)\geqslant 0$. 
Moreover, since $\lam^{\ve}_1(x)<0$ we have 
$$ 
\min_{v\in T_x D}\frac{\imath\pa\db\f(v)}{\curv(E)(v)}<0\,, 
\quad  
\min_{v\in T_x D}\frac{\imath\pa\db\f(v)}{\curv(E)(v)}= 
-\max_{v\in T_x D}\frac{-\imath\pa\db\f(v)}{\curv(E)(v)}\,. 
$$ 
Hence  
\begin{equation}\label{lam1} 
\lam^{\ve}_1\geqslant 1-\ve\chi(\f)M_E(-\f)\quad \text{on $S(\ve)$}\,. 
\end{equation} 
The inequality \eqref{lam1} gives information about the size of $S(\ve)$. 
Indeed, $\lam^{\ve}_1<0$ and \eqref{lam1} entail $\f>-\sqrt{\ve M_E(-\f)}$. 
Thus the integration set is contained in a `corona' of size $\sqrt{\ve\:}$\,: 
\begin{equation}\label{Sepsilon} 
S(\ve)\subset D(\ve/2)\,{\bigcap}\left\lbrace x\in D\,:\,\f(x)>-\sqrt{\ve M_E(-\f)} 
\right\rbrace\,. 
\end{equation} 
Since $\ve\chi(\f)<2\mome$ on $D(\ve/2)$ (see \eqref{not}) we deduce the final  
estimate for the first eigenvalue: 
\begin{equation}\label{lam1fin} 
\abs{\lam^{\ve}_1}\leqslant 2M_E(-\f)\mome-1=:C_2\quad \text{on $S(\ve)$}\,. 
\end{equation} 
We examine now the eigenvalues $\lam^{\ve}_j$ for $j=2,\dotsc,n-1$. 
The min--max principle yields: 
$$ 
\lam^{\ve}_j\leqslant 1+\ve\chi(\f)M_E(\f)+\min_{\substack{F\subset T_x D\\ 
\dim F=j}}\max_{v\in F} 
\frac{\imath\ve\chi^{\prime}(\f)\pa\f\wedge\db\f(v)} 
{\curv(E)(v)}\,. 
$$ 
The minimum in the last expression is $0$ and is attained on some space  
$F\subset\ker\pa\f$. Therefore we get: 
\begin{equation}\label{lam2} 
\abs{\lam^{\ve}_j}\leqslant 1+2\mome M_E(\f)=:C_3\quad 
\text{on $S(\ve)$ for $j=2,\dotsc,n-1$}\,. 
\end{equation} 
The highest eigenvalue satisfies the estimate: 
$$ 
\lam^{\ve}_n\leqslant 1+\ve\chi(\f)M_E(\f)+\ve\chi^{\prime}(\f) 
\max_{v\in T_x D}\frac{\imath\pa\f\wedge\db\f(v)} 
{\curv(E)(v)}\,. 
$$ 
The inequalities: $\ve\chi(\f)<2\mome$ 
and $\ve\chi^{\prime}(\f)\leqslant(2\mome)^{3/2}\ve^{-1/2}$  
hold on $D(\ve/2)$ 
(the last one since $\chi^{\prime}(\f)=-\f^{-3}$). 
We introduce the short notation: 
$$ 
M^{\ve}_{E}(\pa\f)=M_{\ce}(\imath\pa\f\wedge\db\f,K_\ve)\,, 
$$ 
where $K_\ve:=\overline{D}\smallsetminus 
\left\lbrace x\in D\,:\,\f(x)>-\sqrt{\ve M_E(-\f)} 
\right\rbrace$. It is clear that $M^{\ve}_{E}(\pa\f)$ converges to 
$M_{E}(\pa\f)$ for $\ve\longrightarrow 0$. With this notation, 
\begin{equation}\label{lamn} 
\abs{\lam^{\ve}_n}\leqslant 1+2\mome M_E(\f)+ 
\ve^{-1/2}(2\mome)^{3/2}M^{\ve}_{E}(\pa\f)\quad \text{on $S(\ve)$}\,. 
\end{equation} 
At this point we may return to $\abs{I_\ve}$ and use the obvious inequality 
$$ 
\abs{I_\ve}\leqslant 
(2\pi)^{-n}\operatorname{Vol}_{\mathbf{c}(E)}(S(\ve))\, 
\sup_{S(\ve)}\abs{\lam^{\ve}_1}\,\abs{\lam^{\ve}_2} 
\dotsm\abs{\lam^{\ve}_n} 
$$ 
where $\operatorname{Vol}_{\mathbf{c}(E)}$ represents the volume with  
respect to the metric $\ce$. We need to find a bound only for the volume. 
Taking into account \eqref{Sepsilon}, 
\begin{multline}\label{vol} 
\operatorname{Vol}_{\mathbf{c}(E)}(S(\ve))\leqslant 
\sqrt{\ve\,}\left(\sqrt{M_E(-\f)}-\sqrt{(2\mome)^{-1}}\right)\times\\ 
\times\sup\left\lbrace\int_{\{\f=c\}}\frac{dS_E}{\abs{d\f}_E}\,:\, 
c\in\left[-\sqrt{\ve M_E(-\f)},-\sqrt{\ve(2\mome)^{-1}}\right]\right\rbrace 
\end{multline} 
Relations \eqref{lamn} and \eqref{vol} yield: 
\begin{equation*} 
\begin{split} 
\limsup_{\ve\longrightarrow 0}\:&\operatorname{Vol}_{\mathbf{c}(E)}(S(\ve)) 
\sup_{S(\ve)}\abs{\lam^{\ve}_n}\\ 
&\leqslant 
\left(\sqrt{2\mome M_E(-\f)}-1\right)\,2\mome M_E(\pa\f) 
\int_{\pa D}\frac{dS_E}{\abs{d\f}_E}\\ 
&=C_1\,C_4\int_{\pa D}\frac{dS_E}{\abs{d\f}_E} 
\end{split} 
\end{equation*} 
Using \eqref{lam1fin} and \eqref{lam2} we conclude 
\begin{equation}\label{Iepsilon} 
\limsup_{\ve\longrightarrow 0}\;\abs{I_\ve}\leqslant 
(2\pi)^{-n}\,C_1\,C_2\,C^{n-2}_3\,C_4\int_{\pa D}\frac{dS_E}{\abs{d\f}_E} 
\end{equation} 
We are ready to let $\ve\longrightarrow 0$ in \eqref{morse2} and we use 
%the rules $\limsup_{\ve\longrightarrow 0}(a_\ve+b_\ve)\geqslant 
%\limsup_{\ve\longrightarrow 0}a_\ve+\liminf_{\ve\longrightarrow 0}b_\ve$ 
%and $\liminf_{\ve\longrightarrow 0}(-c_\ve)=-\limsup_{\ve\longrightarrow 0}c_\ve$. 
%We apply now 
\eqref{morse3} and \eqref{Iepsilon}. In \eqref{morse3} we can further let 
the compact $L$ exhaust $D$. This proves \eqref{conc1} and with it the  
Existence Criterion. 

%%%%%%%%%%%%%%%%%%%%%%%%%%%%%%%%%%%%%%%%%%%%%%%%%%%%%%%%%%%%%%%%%%%%%%%%%%%%%%
\section{Perturbation of line bundles}

In this section we discuss the relation between the perturbation of the 
complex structure of a line bundle and the perturbation of the complex  
structure on the base manifold. This requires a glance to the corresponding  
section of Lempert's article \cite{Le}. 
Let us consider a compact complex manifold $Y=(Y, \ci)$ with boundary 
endowed with a complex structure $\ci$.  
Let $Z$ be a smooth divisor in $Y$. Denote as usual by $[Z]$ the  
associated line bundle. 
We are interested in the effect of a small perturbation of $\ci$ on $Y$ 
on the complex structure of $[Z]$ or of the canonical bundle $K_Y$  
\underline{over a compact set $D\Subset Y$}.  
This will suffice for the proof of the Stability Theorem. 
Indeed, denote by $E$ a positive line bundle on a concave manifold $Y$ 
and assume that for a small perturbation $\ci^\prime$ of $\ci$  
there exists a perturbation $E^\prime$ of $E$ such that the curvature 
forms of $E$ and $E^\prime$ are close on a sublevel set $D$.  
Then the right hand--side terms  
in \eqref{conc1} calculated for  $\ci$ and $\ci^\prime$  
are also close. If one is positive so is the other and both manifolds 
$D$ and $D^\prime$ (and therefore $Y$ and $Y^\prime$) are Moishezon. 
 
Let us remark that not every perturbation of the complex structure on $Y$ 
lifts to a perturbation of $[Z]$. We need the hypothesis that the tangent 
space $T(Z)$ is $\ci^\prime$ invariant.  
Then $Z$ is a divisor in the new manifold $Y^\prime=(Y,\ci^\prime)$ 
and we consider the associated bundle $[Z]^\prime$. 
Of course any perturbation of $\ci$ 
lifts to a perturbation of the canonical line bundle. 
 
The next Lemma is a ``small perturbation'' of Lemma 4.1 of Lempert \cite{Le}. 
In the latter a compact divisor $Z\subset\operatorname{Int}Y$ is considered 
whereas in our case we deal with a divisor which may cut the boundary. 
However, since we are interested in the effect of the perturbation just on a 
compact set the proof is the same. 
We use the $\cC^\infty$ topology on the spaces of tensors defined on $Y$ 
and also on spaces of restrictions of tensors to compact subsets of $Y$. 
We say that two tensors are close when they are close in the $\cC^\infty$ 
topology. 
 
\begin{lemma} 
Let $(Y,\ci)$ be a compact complex manifold, $Z$ a smooth divisor in $Y$  
and $D\Subset Y$. There exists a finite covering 
$\cu=\{U_\alpha\}_{\alpha\in A}$ of $D$ and a multiplicative cocycle  
$\{g_{\alpha\beta}\in\co_\ci 
(\,\overline{U}_\alpha\cap\overline{U}_\beta\,)\,:\,\alpha,\beta\in A\}$ defining the  
bundle $E=[Z]$ in the vicinity of $D$, with the following property. 
If $\ci^\prime$ is another complex structure on $Y$ close to $\ci$ 
such that $T(Z)$ rests $\ci^\prime$ invariant, 
the bundle $E^\prime$ determined by $Z$ in the structure $\ci^\prime$ 
can be defined in the vicinity of $D$ by the cocycle  
$\{g^\prime_{\alpha\beta}\in\co_{\ci^\prime}(\,\overline{U}_\alpha\cap\overline{U} 
_\beta\,)\,:\,\alpha,\beta\in A\}$  
such that $g^\prime_{\alpha\beta}$ will be as close as we please 
to $g_{\alpha\,\beta}$ on $\overline{U}_\alpha\cap\overline{U}_\beta$  
assuming $\ci^\prime$ and $\ci$ are sufficiently close. 
\end{lemma} 
\begin{proof} 
We remind for the sake of completeness the construction of the cocycles. 
For every point of $Y\cap\overline{D}$ there exists an open neighbourhood  
$U$ in $Y$ and a $\ci$--biholomorphism $\psi_U$ of some neighbourhood of  
$\overline{U}$ into $\C^n$, $n=\dim Y$, such that $\psi_{U}(U)$ is the unit 
polydisc and $\psi_{U}(Z)\subset\{z\in\C^n\,:\,z_1=0\}$.  
Let $\{U_\alpha\}_{1\leqslant\alpha\leqslant m}$ be a finite covering consisting of sets 
$U$ as above and for each $\alpha$ denote by $\psi_\alpha$ the corresponding 
biholomorphism.  
We select further an open set $U_0\Subset Y\smallsetminus Z$ such that  
$\cu=\{U_\alpha\}_{0\leqslant\alpha\leqslant m}$ is a covering of $\overline{D}$. 
For every $1\leqslant\alpha\leqslant m$ we select a smooth strictly  
pseudoconvex Stein domain $U^{\ast}_\alpha\supset U_\alpha$ such that  
$\psi_\alpha$ is biholomorphic in the neighbourhood of $U^{\ast}_\alpha$. 
Set moreover $U^{\ast}_0=U_0$. 
We construct a cocycle defining $E=[Z]$ in the open set $\cup_{\alpha}U^\ast_{\alpha}$ 
as follows. First define functions $g_\alpha$ such that $g_0$ is identically 
$1$ on $U_0$ and $g_\alpha=z_1\circ\psi_\alpha$ for $\alpha\geqslant 1$. 
The bundle $E$ is defined in the vicinity of $D$ by the $\ci$ holomorphic 
multiplicative cocycle $\{g_{\alpha\beta}\}$ where $g_{\alpha\beta}=g_{\alpha}/g_{\beta}$. 
Note that $g_{\alpha\beta}$ is holomorphic on a neighbourhood of  
$\overline{U}^{\ast}_\alpha\cap\overline{U}^{\ast}_\beta\supset 
\overline{U}_\alpha\cap\overline{U}_\alpha$. 
 
Let $\ci^\prime$ be a complex structure as in the statement. Then $Z$ is a complex 
hypersurface in the new structure and defines a line bundle $E^\prime$. 
We describe next the cocycle of $E^\prime$. 
The hypothesis on the sets $U^{\ast}_\alpha$ allows the use of a theorem of  
Hamilton \cite{Ha} for $U^{\ast}_\alpha$. The theorem asserts that for a small 
perturbation $\ci^\prime$ of the complex structure on a neighbourhood 
of $\overline{U}^{\ast}_\alpha$ there is a $\ci^\prime$ biholomorphism  
$\psi^\prime_\alpha$ of a neighbourhood of $\overline{U}^{\ast}_\alpha$ 
into $\C^n$ close to $\psi_\alpha$. As shown in \cite{Le} we can even assume 
$\psi^\prime_\alpha(Z)\subset\{z\in\C^n\,:\,z_1=0\}$. 
Set $g^\prime_0$ to be identically 
$1$ on $U_0$ and $g^\prime_\alpha=z_1\circ\psi^\prime_\alpha$ for $\alpha\geqslant 1$. 
Then put $g^\prime_{\alpha\beta}=g^\prime_\alpha/g^\prime_\beta$. 
Since $\psi_\alpha$ and $\psi^\prime_\alpha$ are close, $g^\prime_\alpha$ is 
$\ci^\prime$ holomorphic on a neighbourhood of $\overline{U}^{\ast}_\alpha$ 
and $g^\prime_{\alpha\beta}$ is $\ci^\prime$ holomorphic on a neighbourhood  
of $\overline{U}^{\ast}_\alpha\cap\overline{U}^{\ast}_\beta$.   
The cocycle $\{g^\prime_{\alpha\beta}\}$ defines $E^\prime$ in the open set 
$\cup_{\alpha}U^\ast_{\alpha}$. 
 
The functions $g_\alpha$ and $g^\prime_\alpha$ are close on  
$\overline{U}_\alpha$. We can now repeat the arguments 
from \cite{Le} to show that $g_{\alpha\beta}$ and $g^\prime_{\alpha\beta}$ 
are also close on $\overline{U}_\alpha\cap\overline{U}_\beta$. 
\end{proof} 
\begin{lemma}\label{curbura} 
Let $(Y,\ci)$, $Z$ and $D\Subset Y$ be as in the preceding Lemma. 
Assume that $[Z]$ is endowed with a hermitian metric $h$.  
If $\ci^\prime$ is another complex structure on $Y$ close to $\ci$, 
leaving $T(Z)$ invariant, there exists a hermitian metric $h^\prime$ on  
the line bundle $[Z]^\prime$ near $D$  
such that the curvature form $\curv([Z]^{\prime})$ 
will be as close as we please to $\curv([Z])$ on $D$ 
assuming $\ci^\prime$ and $\ci$ are sufficiently close. 
\end{lemma} 
\begin{proof} 
We can define a smooth bundle isomorphism $[Z]\longrightarrow [Z]^{\prime}$  
in the vicinity of $D$ 
by resolving the smooth 
additive cocycle $\log(g^\prime_{\alpha\beta}/g_{\alpha\beta})$  
in order to find smooth functions $f_\alpha$, close to $1$  
on a neighbourhood of ${\overline U}_\alpha$ such that  
$g^\prime_{\alpha\beta}=f_\alpha\,g_{\alpha\beta}\,f^{-1}_\beta$.
Then the isomorphism between $[Z]$ and $[Z]^\prime$ is defined by  
$f=\{f_\alpha\}$. 
The metric $h$ is given in terms of the covering $\cu$ by a collection  
$h=\{h_\alpha\}$ of smooth
strictly positive functions satisfying the relation  
$h_\beta=h_\alpha\,\abs{g_{\alpha\beta}}$. We define a hermitian metric
$h^\prime=\{h^\prime_\alpha\}$ on $[Z]^\prime$ by $h^\prime_\alpha= 
h_\alpha\,\abs{f^{-1}_\alpha}$; $h^\prime_\alpha$ is close to $h_\alpha$ 
on $D$.
The curvature form of $[Z]^\prime$ has the form 
$$\frac{\imath}{2\pi}{\bf c}([Z]^\prime)= 
\frac{1}{4\pi}\,d\circ{\ci^\prime}\circ d\,(\log h^\prime_\alpha)\,.$$
Therefore, when $\ci^\prime$ is sufficiently close to $\ci$, 
$\frac{\imath}{2\pi}{\bf c}([Z]^\prime)$ 
is close to $\frac{\imath}{2\pi}{\bf c}([Z])$ on $D$. 
\end{proof} 
In the same vein we study the perturbation of the canonical bundle. 
\begin{lemma}\label{curbura1} 
Let $(Y,\ci)$ and $D\Subset Y$ be as above. 
Assume $K_Y$ is endowed with a hermitian metric $h$.  
If $\ci^\prime$ is another complex structure on $Y$ close to $\ci$, 
there exists a hermitian metric $h^\prime$ on $K_{Y^\prime}$ near $D$ 
such that the curvature form $\curv(K_{Y^\prime})$ 
will be as close as we please to $\curv(K_Y)$ on $D$ 
assuming $\ci^\prime$ and $\ci$ are sufficiently close. 
\end{lemma} 
\begin{proof} 
We find as before a finite covering $\cu=\{U_\alpha\}_{\alpha\in A}$ of 
$\overline{D}$ and biholomorphisms $\psi_\alpha$ defined in a neighbourhood 
of $\overline{U}_\alpha$ which map $U_\alpha$ onto the unit polydisc in $\C^n$.  
For every $\alpha\in A$ we select a smooth strictly  
pseudoconvex Stein domain $U^{\ast}_\alpha\supset U_\alpha$ such that  
$\psi_\alpha$ is biholomorphic in the neighbourhood of $U^{\ast}_\alpha$. 
The canonical bundle $K_Y$ is defined in the vicinity of $D$ by  
$g_{\alpha\beta}=\det\left(\pa\psi_\alpha/\pa\psi_\beta\right) 
=\det\left(\pa\left(\psi_\alpha\circ\psi^{-1}_\beta\right)/\pa w\right)$ 
which is $\ci$--holomorphic on a neighbourhood of  
$\overline{U}^{\ast}_\alpha\cap\overline{U}^{\ast}_\beta\supset 
\overline{U}_\alpha\cap\overline{U}_\alpha$.  
Here $w$ are the canonical coordinates on $\C^n$. 
We apply as before Hamilton's theorem and obtain $\ci^\prime$ 
biholomorphisms $\psi^\prime_\alpha$ in a neighbourhood of 
$\overline{U}^{\ast}_\alpha$ into $\C^n$ close to $\psi_\alpha$. 
 
The cononical bundle $K_{Y^\prime}$ is defined in the vicinity of 
$D$ by 
$g^\prime_{\alpha\beta}=\det\left(\pa\psi^{\prime}_{\alpha}/ 
\pa\psi^{\prime}_\beta\right)$. 
Since $\psi^\prime_\alpha$ is close to $\psi_\alpha$ we see that  
$g^\prime_{\alpha\beta}$ is close to $g_{\alpha\beta}$ on  
$\overline{U}_\alpha\cap\overline{U}_\alpha$. 
By repeating the arguments in the proof of Lemma \ref{curbura} we conclude. 
\end{proof} 
 
%%%%%%%%%%%%%%%%%%%%%%%%%%%%%%%%%%%%%%%%%%%%%%%%%%%%%%%%%%%%%%%%%%%%%%%%%%%%%%%%%
\section{The Stability Theorem}

In this section we prove the Stability Theorem. Let us consider a 
compact manifold $M$, $\dim M=n$, and a complex submanifold $A$ of dimension $q$.
Then $X=M\smallsetminus A$ is $(q+1)$--concave. Let us remind the construction 
of an exhaustion function. Select a finite covering $\hat{\cu}=\{U_\alpha\}_
{\alpha\geqslant 1}$ of $A$
with coordinate domains such that if the coordinates in $U_\alpha$ are
$z_\alpha=(z^1_\alpha,z^2_\alpha,\cdots,z^n_\alpha)$ we have
$A\cap U_\alpha=\{z\in U_\alpha\,:\,z^{q+1}_\alpha=\cdots=z^n_\alpha=0\}$.
Set 
%$\f_\alpha(z)=(\,z^{q+1}_\alpha,\cdots,z^n_\alpha)$ and
$\f_{\alpha}(z)=\sum_{q+1}^n\abs{z^{j}_\alpha}^2$. 
Choose a relatively compact open set $U_0\Subset M\smallsetminus A$
such that $\cu=\{U_0\}\cup\hat{\cu}=\{U_\alpha\}_
{\alpha\geqslant 0}$ is a covering of $M$ and set 
$\f_0\equiv 1$ on $U_0$. Let $\{\rho_\alpha\}_{\alpha\geqslant 0}$ 
be a partition of unity subordinated to $\cu$. Define $\f=\f_A
=\sum_{\alpha\geqslant 0}\rho_\alpha\f_\alpha$.
The function $\f$ enjoys the following properties:
\begin{enumerate}
\item $\f\in\cC^\infty(M)$, $A=\{\f=0\}$ and $\f\geqslant 0$.
\item For any $c>0$ we have $\{\f>c\}\Subset M\smallsetminus A$.
\item $\pa\db\f=\sum_\alpha\left(\rho_\alpha\pa\db\f_\alpha+
\f_\alpha\pa\db\rho_\alpha+\pa\rho_\alpha\wedge\db\f_\alpha+
\pa\f_\alpha\wedge\db\rho_\alpha\right)$ where\\
$\pa\db\f_\alpha=2\sum_{q+1}^n dz^{j}_\alpha\wedge d\overline{z}^j_\alpha$.
\end{enumerate}
For $z\in A$, $\pa\db\f(z)=\sum_{\alpha}\rho_{\alpha}(z)
\pa\db\f_{\alpha}(z)$
has $n-q$ positive eigenvalues. Hence $\pa\db\f$ has $n-q$ positive 
eigenvalues in a neighbourhood of $A$. Moreover $\pa\db\f$ is positive
semidefinite on $A$.
Let us construct a hermitian metric on $M$ which is ``small'' in the 
normal direction to $A$ (near $A$) and ``large'' in the tangential direction
to $A$. We can consider on each $U_\alpha$ the metric
$\delta^{-1}\sum_1^q dz^{j}_\alpha\wedge d\overline{z}^j_\alpha+
\delta\sum_{q+1}^ndz^{j}_\alpha\wedge d\overline{z}^j_\alpha$, ($\delta>0$),
and then patch these metrics together with the partition of unity to obtain a metric $\om_\delta$
on $M$. 
Let $\gamma^{\delta}_1\leqslant\gamma^{\delta}_2\leqslant\dotsb\leqslant\gamma^{\delta}_n$
be the eigenvalues of $\imath\pa\db\f$ with respect to $\om_\delta$.
For $\delta$ sufficiently small there exists a neighbourhood $U_\delta$
of $A$ such that on $U_\delta$, $\gamma^{\delta}_j\geqslant-O(\delta)$ for
$j=1,\dotsc,q$ and $\gamma^{\delta}_j\geqslant O(\delta^{-1})$ for
$j=q+1,\dotsc,n$. Therefore we can choose $\delta$ such that
on $U_\delta$, $\gamma^{\delta}_j\geqslant-1$ for
$j=1,\dotsc,q$ and $\gamma^{\delta}_j\geqslant 2n-3$ for
$j=q+1,\dotsc,n$.

Let us consider now the domains $X_c=\{\f>c\}$ for $c>0$ sufficiently small.
If $\codim A\geqslant 3$ the domains $X_c$ admit as definition function
$c-\f$ whose complex hessian has $3$  negative eigenvalues in the vicinity 
of $\pa X_c$. If $M$ possesses a positive line bundle we are in the 
conditions of the Existence Criterion.
Note that the metric $\om_\delta$ satisfies Property \ref{metric}
for all $X_c$ with $c$ sufficiently small. 
For technical reasons we construct a metric $\om$ as follows. Consider
the real part $g_\delta$ of the hermitian metric $\om_\delta$\,. Thus
$g_\delta$ is a riemannian metric on $M$. Take a hermitian metric
$\om$ whose real part $g$ satisfies $g(u,v)=g_{\delta}(u,v)+
g_{\delta}(\ci u,\ci v)$ ($u,v\in\C\otimes T(M)$)
where $\ci$ is the complex structure of $M$.
If $\delta$ is sufficiently small $\om$ still satisfies Property 
\ref{metric}. From now on we fix such a metric
$\om$ on $M$. The constants $\mome$ are calculated with respect to this 
metric. 
\begin{lemma}
Assume that $M$ is a projective manifold and $E$ is a positive
line bundle over $M$. Let $A$ be a submanifold with $\codim A\geqslant3$.
Then for sufficiently small regular values $c>0$ we have
\begin{equation}\label{integral}
\int_{X_c}\left(\tfrac{\imath}{2\pi}\mathbf{c}(E)\right)^n\,
>\,C(c-\f,E)\int_{\pa X_c}\frac{dS_{E}}{|d\f|_E}
\end{equation}
where $C(c-\f,E)$ has been introduced in \eqref{const}.
\end{lemma}
\begin{proof}
Remark first that the constant $C(c-\f,E)$ converges to $0$ for 
$c\longrightarrow 0$. Indeed, $\pa\db(c-\f)=-\pa\db\f$ so the constants
$M_E(c-\f)$, $M_E(\f-c)$ and $\mome$ are bounded for $c$ running in
a compact interval since $\pa\db\f$ and $E$ are defined over all $M$.
We observe further that $d\f(z)\longrightarrow 0$ when 
$z\longrightarrow  A$ (in fact $d\f\restriction_A=0$). Hence
$M_E(\pa(c-\f)\wedge\db(c-\f),\pa X_c)$ converges to $0$ 
(and with it $C(c-\f)$) when $c$ goes to $0$. 
Examine now the term
$$\int_{\pa X_c}\frac{dS_E}{\abs{d\f}_E}\:.$$
Although $\abs{d\f_E}\longrightarrow 0$ for $z\longrightarrow A$
this integral goes to $0$ too for $c\longrightarrow 0$.
Indeed, since $A$ has codimension $\geqslant 3$ we have
$$\int_{\pa X_c} dS_E=\int_{\{\f=c\}} dS_E=O(c^5)\,,\quad c\longrightarrow 0\,.$$
On the other hand for a regular value $c$ of $\f$,
$$\left\vert d\f\restriction_{\pa X_c}\right\vert
=O(c)\,,\quad c\longrightarrow 0\,.$$
We infer
$$\int_{\pa X_c}\frac{dS_E}{\abs{d\f}_E}=O(c^4)\,,\quad c\longrightarrow 0\,.$$
for regular values $c$ of $\f$.
In conclusion the boundary integral in \eqref{integral} goes to $0$
as $c\longrightarrow 0$. The domain integral in \eqref{integral}
being bounded from below by a positive constant the Lemma follows. 
\end{proof}

At this stage we can prove the Stability Theorem. Let us consider a
smooth domain $Y:=X_c$ for $c$ small enough such that condition
\eqref{integral} holds. Let $\ci^\prime$ be a new complex structure
on $Y$ which leaves $T(Z)$ invariant, for an ample smooth divisor $Z$
on $M$. We apply Lemma \ref{curbura} for the manifold $Y$ and 
a smooth relatively compact set $\overline{D}$ where $D:=X_d$, $d>c$, such that \eqref{integral}
still holds on $X_d$. By hypothesis the bundle $E$ carries
a hermitian metric with positive curvature. Lemma \ref{curbura} shows
that there exists a hermitian metric $h^\prime$ on the
bundle $E^\prime$ near $D$ such that $\curv(E)$ and
$\curv(E^\prime)$ are as close as we please in the $\cC^\infty$
topology on $\overline{D}$ if $\ci$ and $\ci^\prime$ are sufficiently close.
In particular $\curv(E^\prime)$ is positive near $\overline{D}$.
Note that a defining function for $D^\prime$ is still $d-\f$
and its complex hessian will have $3$ negative eigenvalues in the
vicinity of $\pa D^\prime$
for a small perturbation of the complex structure. 

Thus we can apply the Existence Criterion for $D^\prime$ and $E^\prime$.
In order to calculate the constant $C(d-\f,E^\prime)$ we construct first
a metric $\om^\prime$ on $Y$ in the following way.
The metric $\om$ determines a riemannian metric $g$ on $Y$ which was
chosen such that $g(u,v)=g_{\delta}(u,v)+
g_{\delta}(\ci u,\ci v)$ for $u,v\in\C\otimes T(M)$.
We consider then a hermitian metric $\om^\prime$ on $Y^\prime$
with real part
$g^\prime$ where $g^{\prime}(u,v)=g_{\delta}(u,v)+
g_{\delta}(\ci^{\prime}u,\ci^{\prime}v)$
for $u,v\in\C\otimes T(M)$.
The metric $\om^\prime$ satisfies the Property \ref{metric}
with respect to the defining function $d-\f$ of $D^\prime$, provided 
$\ci$ and $\ci^\prime$ are sufficiently close.
Therefore the constants $M_{E^\prime}(d-\f)$, $M_{E^\prime}(\f-d)$,
$M_{\om^\prime}(E^\prime)$ and $M_{E^\prime}(\pa(d-\f),\pa D^\prime)$
are close to the corresponding constants  $M_{E}(d-\f)$, $M_{E}(\f-d)$,
$M_{\om}(E)$ and $M_{E}(\pa(d-\f),\pa D)$ respectively.
This entails that $C(d-\f,E^\prime)$ is close to $C(d-\f,E)$.

It is also clear that $\int_{D^\prime}\left(\tfrac{\imath}{2\pi}
\curv(E^\prime)\right)^n$  
and $\int_{\pa{D^\prime}} dS_{E^\prime}/\abs{d\f}_{E^\prime}$,
are close to the corresponding
integrals on $D$ and $\pa D$ of $\tfrac{\imath}{2\pi}\ce$ and 
$dS_{E^\prime}/\abs{d\f}_{E^\prime}$.
Therefore
\begin{equation}
\tag{\ref{integral}$^\prime$}\label{intprime}
\int_{D^\prime}\left(\tfrac{\imath}{2\pi}\mathbf{c}(E^\prime)\right)^n\,
>\,C(d-\f,E^\prime)\int_{\pa D^\prime}\frac{dS_{E^\prime}}{|d\f|_{E^\prime}}
\end{equation}
By the Existence Criterion
\begin{equation}
\dim H^0(D^\prime,{E^\prime}^k)\gtrsim k^n
\end{equation}
for large $k$ and thus $D^\prime$ and so $Y^\prime$ are Moishezon,
provided $\ci$ and $\ci^\prime$ are sufficiently close.
An entirely analogous argument takes care of the case of perturbation of the
canonical bundle $K_Y$.
This proves the Stability Theorem.

%%%%%%%%%%%%%%%%%%%%%%%%%%%%%%%%%%%%%%%%%%%%%%%%%%%%%%%%%%%%%%%%%%%%%%%%%%%%%%% 

\end{document}